%% file: BergLanz_MTNS_-_revised_arxiv.tex
\documentclass{ifacconf}

\usepackage{amsmath}
\usepackage{amssymb}
\usepackage{natbib}        

\usepackage{pbox}
\usepackage{graphicx}      
\usepackage{tikz}
\usetikzlibrary{arrows}
\usepackage{caption}

\makeatletter
\def\endfigure{\end@float}
\def\endtable{\end@float}
\makeatother
\usepackage[caption=false]{subfig}
\usepackage{etoolbox}
\AtBeginEnvironment{figure}{\setlength{\captionwidth}{0.9\linewidth}}
\AtBeginEnvironment{table}{\setlength{\captionwidth}{0.9\linewidth}}
\newcommand{\cF}{\mathcal{F}}

\newcommand{\B}{\mathfrak{B}}

\newcommand{\N}{{\mathbb{N}}}

\newcommand{\R}{{\mathbb{R}}}
\newcommand{\C}{{\mathbb{C}}}

\newcommand{\Gl}{\mathbf{Gl}}

\DeclareMathOperator{\im}{im}

\DeclareMathOperator{\diag}{diag}

\newcommand{\setdef}[2]{\left\{\ #1\ \left|\ \vphantom{#1} #2 \right.\right\}}

\theoremstyle{definition}
\newtheorem{Rem}{Remark}

\newenvironment{smallpmatrix}
{\left(\begin{smallmatrix}}
{\end{smallmatrix}\right)}

\newenvironment{smallbmatrix}
{\left[\begin{smallmatrix}}
{\end{smallmatrix}\right]}

\begin{document}
\begin{frontmatter}

\title{Output tracking for a non-minimum phase robotic manipulator\thanksref{footnoteinfo}}

\thanks[footnoteinfo]{This work was supported by the German Research Foundation (Deutsche Forschungsgemeinschaft) via the grant BE 6263/1-1.}

\author[First]{Thomas Berger}
\author[First]{Lukas Lanza}

\address[First]{Institut f\"ur Mathematik, Universit\"at Paderborn, Warburger Str.~100, 33098~Paderborn, Germany
\newline (e-mail: \{thomas.berger, lanza\}@math.upb.de).}

\begin{abstract}                
We exploit a recently developed funnel control methodology for linear non-minimum phase systems to design an output error feedback controller for a nonlinear robotic manipulator, which is not minimum phase.
We illustrate the novel control design by a numerical case study, where we simulate end-effector output tracking of the robotic manipulator.
\end{abstract}


\end{frontmatter}

\section{Introduction}
Throughout the last two decades funnel control (introduced in~\cite{IlchRyan02b}) turned out to be a powerful tool for the treatment of tracking problems in various applications, such as speed control of wind turbine systems \cite{Hack14,Hack15a}, termination of fibrillation processes~\cite{BergBrei21}, control of peak inspiratory pressure \cite{PompWeye15}, temperature control of chemical reactor models \cite{IlchTren04}, current and voltage control of electrical circuits \cite{BergReis14a}, adaptive cruise control \cite{BergRaue18,BergRaue20} and control of industrial servo-systems \cite{Hack17} and underactuated multibody systems \cite{BergOtto19}.
First results for robotic manipulators have been derived in~\cite{HackEndi08b, HackKenn12}.
Many of the aforementioned applications are such that their dynamics are minimum phase. Concerning this property, there are some nuances in the literature, see e.g. \cite{IlchWirt13}; we call a system minimum phase, if its internal dynamics (in the linear case the zero dynamics) are bounded-input, bounded-output stable.

\begin{figure}[h!tp]
\begin{center}
\begin{tikzpicture}[scale=0.45]
\tikzset{>=latex}
  \filldraw[color=gray!25] plot[smooth] coordinates {(0.15,4.7)(0.7,2.9)(4,0.4)(6,1.5)(9.5,0.4)(10,0.333)(10.01,0.331)(10.041,0.3) (10.041,-0.3)(10.01,-0.331)(10,-0.333)(9.5,-0.4)(6,-1.5)(4,-0.4)(0.7,-2.9)(0.15,-4.7)};
  \draw[thick] plot[smooth] coordinates {(0.15,4.7)(0.7,2.9)(4,0.4)(6,1.5)(9.5,0.4)(10,0.333)(10.01,0.331)(10.041,0.3)};
  \draw[thick] plot[smooth] coordinates {(10.041,-0.3)(10.01,-0.331)(10,-0.333)(9.5,-0.4)(6,-1.5)(4,-0.4)(0.7,-2.9)(0.15,-4.7)};
  \draw[thick,fill=lightgray] (0,0) ellipse (0.4 and 5);
  \draw[thick] (0,0) ellipse (0.1 and 0.333);
  \draw[thick,fill=gray!25] (10.041,0) ellipse (0.1 and 0.333);
  \draw[thick] plot[smooth] coordinates {(0,2)(2,1.1)(4,-0.1)(6,-0.7)(9,0.25)(10,0.15)};
  \draw[thick,->] (-2,0)--(12,0) node[right,above]{\normalsize$t$};
  \draw[thick,dashed](0,0.333)--(10,0.333);
  \draw[thick,dashed](0,-0.333)--(10,-0.333);
  \node [black] at (0,2) {\textbullet};
  \draw[->,thick](4,-3)node[right]{\normalsize$\lambda$}--(2.5,-0.4);
  \draw[->,thick](3,3)node[right]{\normalsize$(0,e(0))$}--(0.07,2.07);
  \draw[->,thick](9,3)node[right]{\normalsize$\varphi(t)^{-1}$}--(7,1.4);
\end{tikzpicture}
\end{center}
\caption{Error evolution in a funnel $\cF_{\varphi}$ with boundary $\varphi(t)^{-1}$ for $t > 0$.}
\label{Fig:funnel}
\end{figure}
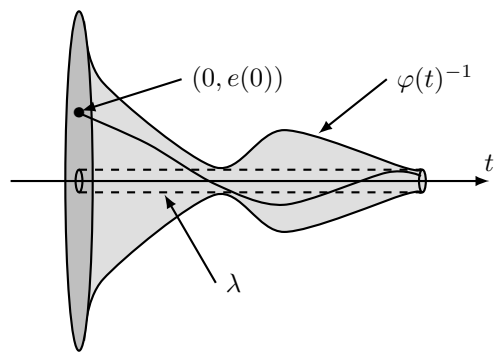

The objective of funnel control is to design a feedback control law such that in the closed-loop system the tracking error~${e(t) = y(t) - y_{\rm ref}(t)}$ evolves within the boundaries of a prescribed performance funnel
\begin{align*}
\mathcal{F}_{\varphi} := \setdef{ (t,e) \in \R_{\geq 0} \times \R^m }{ \varphi(t) \| e\| < 1 },
\end{align*}
which is determined by a so called funnel function~$\varphi$ belonging to a large set of functions
\begin{align*}
\Phi_r := \left\{\varphi \!\in\!  C^r(\R_{\geq 0} \!\to\! \R) \ \vline \
\begin{array}{l}
\varphi, \dot \varphi,\ldots,\varphi^{(r)} \, \text{are bounded,}\\
\varphi (\tau)>0 \, \text{ for all} \, \tau>0, \\
\text{and }  \liminf_{\tau\rightarrow \infty} \varphi(\tau) > 0
\end{array}
\right\},
\end{align*}
where $r\in\N$ is the relative degree of the system. The latter is, roughly speaking the minimal number one has to differentiate the output of a system to obtain the input explicitly; a definition is provided in the subsequent subsection. The boundary of~$\cF_\varphi$ is given by $1/\varphi$, a typical situation is depicted in Fig.~\ref{Fig:funnel}. Note that, since~$\varphi$ is bounded, the performance funnel is bounded away from zero, which means that there exists $\lambda>0$ so that $1/\varphi(t) \ge \lambda$ for all $t>0$.
The requirement of bounded~$\varphi$ can also be waived, see the recent work~\cite{BergIlch21}.

Our long-term aim is to study output tracking for nonlinear non-minimum phase systems with known relative degree via funnel control.
To achieve this, we seek to exploit a recent result for linear non-minimum phase systems in \cite{Berg20a}, where output tracking with prescribed performance is achieved by introducing an auxiliary output. The latter is interpreted as output of a new system with higher relative degree. Using a funnel controller for nonlinear systems with known relative degree as in \cite{BergLe18a} and an appropriate new reference signal, tracking with prescribed performance is achieved.

In the present work, we use the ideas and method developed in~\cite{Berg20a} to design a controller with which we perform a numerical case study
and determine whether the ideas from~\cite{Berg20a} are, in principle, feasible for nonlinear systems. We focus on the presentation of a controller design and its numerical validation; feasibility proofs are a topic of future research.
As a specific example, we consider end-effector output tracking of a robotic manipulator arm. This example is taken from~\cite{SeifBlaj13}. We emphasize that the robotic manipulator is nonlinear and underactuated, i.e., the system has fewer inputs than degrees of freedom. Further, the internal dynamics are not stable, but have a hyperbolic equilibrium. Moreover, if we treat the internal dynamics as a control system, then a flat output does not exist (cf.~\cite{FlieLevi95} for flatness). Therefore, the approach from~\cite{Berg20a} is not directly feasible. To resolve this, the controller design is based on the linearization of the internal dynamics.
Unlike \cite{SeifBlaj13} we use a feedback control law.

A popular
open-loop
 alternative to the approach presented in the present paper is stable system inversion, see~\cite{ChenPade96,DevaChen96}. Based on the given reference trajectory, in this method an open-loop (feedforward) control input is calculated. However, for non-minimum phase systems this requires a reverse-time integration, and hence the computed open-loop control input is non-causal. A different approach is based on so called ideal internal dynamics, see~\cite{GopaHedr93,ShkoShte02}. However, this method requires a trackability assumption and sufficient conditions for its feasibility are not available. Stabilization of non-minimum phase systems by dynamic compensators was considered in~\cite{Isid00} and extended to regulator problems in~\cite{MarcIsid04, NazrKhal09}, where in the latter work extended high-gain observers are used.

\noindent\textbf{Notation}: $\mathbb{N}$ and $\mathbb{R}$ denote the natural and real numbers, resp., $\N_0 = \mathbb{N} \cup \{0\}$ and $\mathbb{R}_{\geq 0} = [0,\infty)$. $C^k(\mathbb{R}_{\geq 0}, \mathbb{R})$ is the (linear) space of $k$-times continuously differentiable functions $f:\R_{\ge 0}\to\R$, and $L^\infty(\mathbb{R}_{\geq 0}, \mathbb{R})$ denotes the (linear) space of Lebesgue-measurable and essentially bounded functions. Moreover, $W^{k,\infty}(\mathbb{R}_{\geq 0},\mathbb{R})$ is the Sobolev space of all functions $f:\R_{\ge 0}\to\R$ with $f,\dot f,\ldots, f^{(k)}\in L^\infty(\mathbb{R}_{\geq 0}, \mathbb{R})$, where $k\in\N$. By $\Gl_n(\R)$ we denote the general linear group of all invertible matrices $A\in\R^{n\times n}$.

\section{System class and Byrnes-Isidori form}

We consider nonlinear systems of the form
\begin{align}
\label{eq:System}
 \dot{x}(t) &= f(x(t)) + g(x(t)) u_d(t), \quad x(0) = x^0 \in \R^n, \nonumber \\
y(t) &= h(x(t)),
\end{align}
with ${f: \R^n \to \R^n}$,  ${g: \R^n \to \R^{n \times m}}$ sufficiently smooth vector fields and ${h: \R^n \to \R^m}$ a sufficiently smooth mapping. Further, we assume that $u_d(t) = u(t) + d(t)$, where~$u: \R_{\geq 0} \to \R^m$ and~$y: \R_{\geq 0} \to \R^m$ denote the input and output, resp., and $d:\R_{\geq 0} \to \R^n$ are bounded disturbances or uncertainties. Disturbances of this kind are expected in real applications, in particular in multibody systems.  Note that the dimensions of input and output coincide.

We assume, that system~\eqref{eq:System} has relative degree~$r\in\N$ in the following sense.
First, recall the definition of the \textit{Lie derivative} of a function~$h$ along a vector field~$f$ at a point~$z \in U \subseteq \R^n$,~$U$ open:
\begin{align*}
\left(L_f h \right)(z) := h'(z) f(z),
\end{align*}
where~$h'$ is the Jacobian of~$h$. We may gradually define~$L_f^k h = L_f(L_f ^{k-1} h)$ with $L_f^0 h = h$. Furthermore, denoting with $g_i(z)$ the columns of $g(z)$ for $i=1,\ldots,m$, we define
\[
    \left(L_g h \right)(z) := [(L_{g_1} h)(z), \ldots,(L_{g_m} h)(z)].
\]
Now, in virtue of~\cite{Isid95}, the system~\eqref{eq:System} is said to have relative degree~$r \in \mathbb{N}$ on~$U$, if for all~$z \in U$ we have:
\begin{align*}
&\forall\, k \in \{ 0,...,r-2 \}: \ (L_g L_f^k h)(z) = 0_{m \times m} \\
\text{and}\quad & \Gamma(z) :=  (L_g L_f^{r-1} h)(z) \in \Gl_m(\R),
\end{align*}
where $\Gamma:U\to\Gl_m(\R)$ denotes the high-frequency gain matrix.

As mentioned above, in the present paper the design of the controller is based on the linearization of the internal dynamics. Therefore, we have to decouple the internal dynamics first. To this end, we transform system~\eqref{eq:System} nonlinearly into Byrnes-Isidori form, for details see e.g.~\cite{Isid95}.
If system~\eqref{eq:System} has relative degree~$r \in \mathbb{N}$ on an open set $U \subseteq \R^n$, then there exists a (local) diffeomorphism~${\Phi: U \to W \subseteq \R^n}$, $W$ open, such that the coordinate transformation $\begin{smallpmatrix} \xi(t) \\ \eta(t) \end{smallpmatrix} = \Phi(x(t))$, $\xi(t)\in\R^{rm}$, $\eta(t)\in\R^{n-rm}$ puts system~\eqref{eq:System} into Byrnes-Isidori form
\begin{align}
\label{eq:BIF}
y(t) &= \xi_1(t), \nonumber \\
\dot{\xi}_1(t) &= \xi_2(t),  \nonumber \\
& \ \ \vdots \\
\dot{\xi}_{r-1}(t) &= \xi_r(t), \nonumber \\
\dot{\xi}_r(t) &= (L_f^r h)\big(\Phi^{-1}\!(\xi(t),\eta(t))\big) \!+\! \Gamma\big(\Phi^{-1}\!(\xi(t),\eta(t))\big) u_d(t), \nonumber \\
\dot{\eta}(t) &= q(\xi(t),\eta(t)) + p(\xi(t),\eta(t))u_d(t). \nonumber
\end{align}
The last equation in~\eqref{eq:BIF} represents the internal dynamics of system~\eqref{eq:System}. The diffeomorphism~$\Phi$ can be represented as
\begin{align}
\label{eq:BIF-trafo}
\Phi(x) =
\begin{smallpmatrix}
h(x) \\
(L_fh)(x) \\
\vdots \\
(L_f^{r-1}h)(x) \\
\phi_1(x) \\
\vdots \\
\phi_{n-rm}(x)
\end{smallpmatrix},
\end{align}
where $\phi_i:U\to\R$, $i=1,\ldots,n-rm$ are such that $\Phi'(z)\in\Gl_n(\R)$ for all $z\in U$. If the distribution~$\im (g(x))$ in~\eqref{eq:System} is involutive, then the functions~$\phi_i$ in~\eqref{eq:BIF-trafo} can additionally be chosen such that
\begin{align}
\label{eq:involutivity}
\forall\, i=1,\ldots,n-rm\ \forall\, z\in U:\ (L_g \phi_i) (z) = 0,
\end{align}
by which $p(\cdot) = 0$ in~\eqref{eq:BIF}, cf.~\cite[Prop.~5.1.2]{Isid95}. Recall from~\cite[Sec.~1.3]{Isid95} that~$\im (g(x))$ is involutive, if for all smooth vector fields $g_1, g_2: \R^n \to \R^n$ with $g_i(x)\in \im (g(x))$ for all $x \in \R^n$ and $i=1,2$ we have that the Lie bracket $[g_1,g_2](x) = g_1'(x) g_2(x) - g_2'(x) g_1(x)$ satisfies $[g_1,g_2](x) \in \im (g(x))$ for all $x\in\R^n$.

\section{Rotational manipulator arm}

\begin{figure}[ht!b]
\hspace{1.5cm}\includegraphics[width=0.7\linewidth]{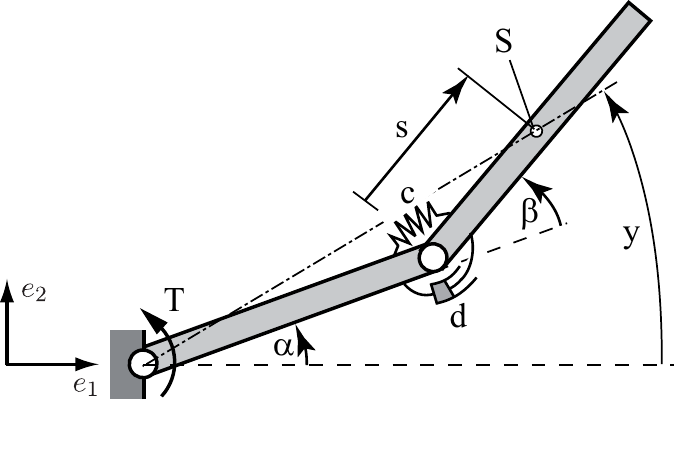}
\caption{Rotational manipulator arm with two links and a passive joint. The figure is taken from~\cite{SeifBlaj13}.}
\label{Fig:Manipulator}
\end{figure}

To illustrate our approach, in the present paper we consider an underactuated rotational manipulator arm as in~\cite{SeifBlaj13}, depicted in Fig. \ref{Fig:Manipulator}.
The manipulator arm consists of two equal links with homogeneous mass distribution,
mass~$m$,
length~$l$ and inertia~${I = l^2 m /12}$.
The two links are coupled via a passive joint consisting of a linear spring-damper combination with spring constant~$c$,
and damping coefficient~$d$.
Using a body fixed coordinate system the tracking point~$S$ on the passive link is described by ${0 \leq s \leq l}$.
The control input is chosen as a torque~$u=T$ acting on the first link.

Define $\B := \setdef{ \beta \in [-\pi,\pi)}{ \cos(\beta) > 2/3 }$ and $U_\beta := \R\times\B\times\R^2$.
Then, as we will see soon, it is reasonable to consider the dynamics of the manipulator on the set~$U_\beta$, i.e., the angle $\beta(t)$ is restricted to~$\B$.
Next, we present the equations of motion of the manipulator arm. With
\begin{align*}
& & M: U_\beta & \to \R^{2 \times 2},\\
&& (x_1,\ldots,x_4) &\mapsto
l^2 m \begin{bmatrix}
\tfrac{5}{3} + \cos(x_2) & \tfrac{1}{3} + \tfrac{1}{2}\cos(x_2) \\ \tfrac{1}{3} + \tfrac{1}{2}\cos(x_2) & \tfrac{1}{3}
\end{bmatrix}, \\
& & f_1 : U_\beta&\to \R, \\
&& (x_1,\ldots,x_4) &\mapsto \tfrac{1}{2} l^2 m x_4 (2 x_3 + x_4) \sin(x_2), \\
& & f_2 : U_\beta &\to \R, \\
&& (x_1,\ldots,x_4) &\mapsto -c x_2 - d x_4 - \tfrac{1}{2} l^2 m x_3^2 \sin(x_2)
\end{align*}
the equations of motion are given by
\begin{align}
\label{eq:motion-manipulator}
&M\big(\alpha(t),\beta(t),\dot\alpha(t),\dot\beta(t)\big) \begin{pmatrix}
\ddot{\alpha}(t) \\ \ddot{\beta}(t)
\end{pmatrix}\nonumber \\
& = \begin{pmatrix} f_1\big(\alpha(t),\beta(t),\dot\alpha(t),\dot\beta(t)\big) \\ f_2\big(\alpha(t),\beta(t),\dot\alpha(t),\dot\beta(t)\big) \end{pmatrix}
 + \begin{bmatrix} 1 \\ 0 \end{bmatrix} u_d(t).
\end{align}
For later use, we compute the inverse of the mass matrix:
\begin{align*}
M(x)^{-1} = \tfrac{36}{(l^2m)^2(16-9\cos^2(x_2))}
\begin{smallbmatrix}
 \tfrac{1}{3} & -\tfrac{1}{3} - \tfrac{1}{2}\cos(x_2) \\ -\tfrac{1}{3} - \tfrac{1}{2}\cos(x_2) & \tfrac{5}{3} + \cos(x_2)
\end{smallbmatrix}
\end{align*}
for $x=(x_1,\ldots,x_4)^\top\in U_\beta$. As an output for~\eqref{eq:motion-manipulator} we consider the auxiliary angle
\begin{equation}
\label{eq:output-manipulator}
y(t) = \alpha(t) + \frac{s}{s+l}\, \beta(t),
\end{equation}
which approximates the position~$S$ on the passive link for a small angle~$\beta$.
Now, in~\cite{SeifBlaj13} it is shown that for this particular output and ${s/l > 2/3}$ the linearized internal dynamics are unstable.
Hence output tracking with ${s=l}$, i.e., end-effector tracking, leads to unstable internal dynamics, which is the case we are interested in.

To find the relative degree of~\eqref{eq:motion-manipulator} we calculate the corresponding Lie derivatives. However, this is not directly feasible, since we require a first-order formulation as in~\eqref{eq:System}. We may rewrite~\eqref{eq:motion-manipulator} in this form with
\begin{align*}
  f : U_\beta&\to \R^4, \\
  x=(x_1,\ldots,x_4) &\mapsto \diag\big(I_2, M(x)^{-1}\big) \begin{smallpmatrix}
    x_3\\ x_4\\ f_1(x)\\ f_2(x) \end{smallpmatrix},\\
    g : U_\beta&\to \R^{4\times 1}, \\
  x=(x_1,\ldots,x_4) &\mapsto \diag\big(I_2, M(x)^{-1}\big) \begin{smallbmatrix}
    0\\ 0\\ 1\\ 0 \end{smallbmatrix},\\
    h : U_\beta&\to \R, \\
  x=(x_1,\ldots,x_4) &\mapsto \left[1, \tfrac{s}{s+l}, 0, 0\right]\, x.\\
\end{align*}
Then, for any $z\in U_\beta$, we obtain the Lie derivatives
\begin{align*}
  (L_g h)(z) &=  \left[1, \tfrac{s}{s+l}, 0, 0\right] \diag\big(I_2, M(z)^{-1}\big) \begin{smallbmatrix}
    0\\ 0\\ 1\\ 0 \end{smallbmatrix} = 0,\\
   (L_g L_f h)(z) &= (L_f h)'(z) g(z) = \big(h' f\big)'(z) g(z)\\
   &= \left( \left[1, \tfrac{s}{s+l}, 0, 0\right] \diag\big(I_2, M(z)^{-1}\big) \begin{smallpmatrix}
    z_3\\ z_4\\ f_1(z)\\ f_2(z) \end{smallpmatrix}\right)' \\
    &\quad \cdot \diag\big(I_2, M(z)^{-1}\big) \begin{smallbmatrix}
    0\\ 0\\ 1\\ 0 \end{smallbmatrix} \\
   & = \left[0, 0, 1, \tfrac{s}{s+l}\right] \diag\big(I_2, M(z)^{-1}\big) \begin{smallbmatrix}
    0\\ 0\\ 1\\ 0 \end{smallbmatrix}\\
   & = \left[1, \tfrac{s}{s+l}\right]  M(z)^{-1} \begin{bmatrix} 1\\ 0 \end{bmatrix}.
\end{align*}
The latter expression is the high-frequency gain matrix
\begin{align*}
\Gamma: U_\beta&\to \Gl_1(\R), \\
x&\mapsto\left[1, \tfrac{s}{s+l}\right]  M(x)^{-1} \begin{bmatrix} 1\\ 0 \end{bmatrix} \\
&\quad\ \, = \tfrac{36}{l^2m(16-9\cos^2(x_2))}\left[ \tfrac{1}{3} - \tfrac{s}{l+s} \left(\tfrac{1}{3} + \tfrac{1}{2}\cos(x_2) \right) \right]
\end{align*}
which is invertible for~$s=l$ since $x_2\in\B$ for all $x\in U_\beta$. Therefore, the relative degree of~\eqref{eq:motion-manipulator},~\eqref{eq:output-manipulator} with $s=l$ is $r=2$ on $U_\beta$.
In passing, we mention the observation that for~$s=l$ and $x\in \R\times [-\pi,\pi)\times\R^2$ we have $\Gamma(x) < 0$ if, and only if, $x\in U_\beta$, i.e., $x_2\in \B$. Since we seek to consider an open area around the equilibrium $x=0$, where the relative degree is well defined, $U_\beta$ is the largest set where this is true.

\section{Output tracking control}

In this section we present the novel controller design with which we perform the numerical case study in Section~\ref{NumericalExample}.
First, as a motivation, we briefly recall the controller from~\cite{Berg20a} for linear non-minimum phase systems.
Then, we derive the representation of system~\eqref{eq:motion-manipulator}, \eqref{eq:output-manipulator} in Byrnes-Isidori form in order to isolate the internal dynamics. Based on the linearization of the internal dynamics we design a controller for end-effector output tracking of the manipulator.

\subsection{Methodology}\label{Ssec:Methodology}

Recently, a controller for linear non-minimum phase systems was developed in~\cite{Berg20a}. There, a differentially flat output~$y_{\rm new}$ for the unstable part of the internal dynamics is introduced. Recall that all state and input variables can be parameterized in terms of a flat output, if it exits, cf.~\cite{FlieLevi95}. Interpreting the new output as output for a system with higher relative degree, the unstable part of the internal dynamics is removed. Then applying a funnel control law as in~\cite{BergLe18a} to the system with appropriate new reference signal leads to tracking with prescribed performance.

In order to adopt this approach for the robotic manipulator~\eqref{eq:motion-manipulator},~\eqref{eq:output-manipulator} we first compute its internal dynamics based on the Byrnes-Isidori form~\eqref{eq:BIF}. We then observe that the internal dynamics depend nonlinearly on~$\dot y$ and, when this is considered as the input, the unstable part does not have a flat output. Hence, we linearize the internal dynamics around the equilibrium and apply the above mentioned controller design based on this linearization. Since derivatives of the new output, which is defined in this way, are required we need to replace the variables obtained via the linearization with the variables of the original system. As an alternative, we present an approach where the derivatives are calculated using a high-gain observer.


\subsection{Internal dynamics of the manipulator}

Henceforth we consider end-effector tracking and set~$s=l$. In order to obtain the internal dynamics, we transform system~\eqref{eq:motion-manipulator} with output~\eqref{eq:output-manipulator} into Byrnes-Isidori form. For $U_\beta\subseteq\R^4$ as above, define $\Phi:U_\beta\to\R^4$ as in~\eqref{eq:BIF-trafo} with $\phi_i:U_\beta \to\R$, $i=1,2$. Since the distribution
\[
    \im g(x) = \{0\}^2 \times \im M(x)^{-1} \begin{bmatrix} 1\\ 0 \end{bmatrix}
\]
is one-dimensional and hence obviously involutive, we may choose~$\phi_1$ and~$\phi_2$ such that
\begin{equation}\label{eq:cond-phi-i}
\forall\, x\in U_\beta:\    0 = (L_g \phi_i)(x) = \phi_i'(x) \diag\big(I_2, M(x)^{-1}\big) \begin{smallbmatrix}
    0\\ 0\\ 1\\ 0 \end{smallbmatrix}.
\end{equation}
Similar to the recent work~\cite{Lanz21},
we investigate the ansatz
\begin{align*}
  \phi_1(x_1,\ldots,x_4)  &= \tilde \phi_1(x_1,x_2),\\
  \text{and}\quad \phi_2(x_1,\ldots,x_4)  &= \tilde \phi_2(x_1,x_2) \begin{pmatrix} x_3\\ x_4\end{pmatrix}
\end{align*}
for $\tilde\phi_1:U_1\to\R$, $\tilde\phi_2: U_1\to\R^{1\times 2}$ and $U_1 = \R\times\B$. Since we require the transformation~$\Phi$ to be a local diffeomorphism, its Jacobian has to be invertible on~$U_\beta$. This is the case if, and only if,
\begin{align}
\label{eq:Invertibility-Condition}
& \forall\, x\in U_\beta:\ \Phi'(x) = \begin{bmatrix}
h'(x) \\ (h' f)'(x)\\ \phi_1'(x)\\ \phi_2'(x)\end{bmatrix} \nonumber\\
&\qquad\qquad\ = \begin{bmatrix} [1, 1/2] & 0\\ 0 & [1, 1/2]\\ \tilde\phi_1'(x_1,x_2) & 0\\ \ast & \tilde \phi_2(x_1,x_2)\end{bmatrix}\in\Gl_4(\R)  \nonumber\\
\iff\ \ & \forall\, q\in U_1:\ \begin{bmatrix}
[1, 1/2] \\
\tilde \phi_1(q)\end{bmatrix}\in\Gl_2(\R)  \nonumber\\
&\qquad\qquad \quad \wedge \quad
\begin{bmatrix}
[1, 1/2] \\
\tilde \phi_2(q)
\end{bmatrix}\in\Gl_2(\R).
\end{align}
In order to satisfy conditions~\eqref{eq:cond-phi-i} and~\eqref{eq:Invertibility-Condition} we choose
\begin{align}
\label{eq:choice-phi}
\tilde \phi_1: U_1 &\to \R,\ q \mapsto q_2, \nonumber \\
\tilde \phi_2: U_1 &\to \R^{1 \times 2},\ q \mapsto  \left[\tfrac{1}{3} + \tfrac{1}{2}\cos(q_2), \tfrac{1}{3}\right].
\end{align}
With this choice, clearly~\eqref{eq:cond-phi-i} is satisfied. In order to verify~\eqref{eq:Invertibility-Condition} we calculate
\begin{align*}
\begin{bmatrix}
[1, 1/2] \\
\tilde \phi_1'(q)
\end{bmatrix} = \begin{bmatrix} 1 & \tfrac{1}{2} \\ 0 & 1 \end{bmatrix} \ \text{and} \
\begin{bmatrix}
[1, 1/2] \\
\tilde \phi_2(q)
\end{bmatrix} = \begin{bmatrix} 1 & \tfrac{1}{2} \\ \tfrac{1}{3} + \tfrac{1}{2}\cos(q_2) & \tfrac{1}{3} \end{bmatrix}
\end{align*}
for $q\in U_1$, where the latter is invertible since $q_2 \in \B$. We may now infer that~\eqref{eq:motion-manipulator}, \eqref{eq:output-manipulator} can be transformed into Byrnes-Isidori form with the particular choice for~$\tilde{\phi}_1, \tilde{\phi}_2$ as in~\eqref{eq:choice-phi}.
Substituting the respective expressions via
\begin{equation}\label{eq:BIF-Manipulator}
\begin{pmatrix} y(t)\\ \dot y(t)\\ \eta_1(t)\\ \eta_2(t)\end{pmatrix} = \Phi\begin{pmatrix} \alpha(t)\\ \beta(t)\\ \dot\alpha(t)\\ \dot\beta(t)\end{pmatrix}
\end{equation}
and rearranging yields
\begin{align*}
\begin{pmatrix} y(t) \\ \eta_1(t) \end{pmatrix} &= \begin{bmatrix} 1 & \tfrac{1}{2} \\ 0 & 1 \end{bmatrix} \begin{pmatrix} \alpha(t) \\ \beta(t) \end{pmatrix}, \\
\begin{pmatrix} \dot{y}(t) \\ \eta_2(t) \end{pmatrix} &= \begin{bmatrix} 1 & \tfrac{1}{2} \\ \tfrac{1}{3} + \tfrac{1}{2}\cos(\beta(t)) & \tfrac{1}{3} \end{bmatrix} \begin{pmatrix} \dot{\alpha}(t) \\ \dot{\beta}(t) \end{pmatrix}.
\end{align*}
Solving the above equations for $y, \eta_1$ and $\dot y, \eta_2$, we may now formulate the internal dynamics as in~\eqref{eq:BIF} as follows:
\begin{align}
\label{eq:ID}
\dot{\eta}_1(t) &= g_{1,0}(\eta_1(t),\eta_2(t)) + g_{1,1}(\eta_1(t),\eta_2(t))\, \dot{y}(t), \nonumber\\
 \dot{\eta}_2(t)&= g_{2,0}(\eta_1(t),\eta_2(t)) + g_{2,1}(\eta_1(t),\eta_2(t))\, \dot{y}(t)\\
 &\quad  + g_{2,2}(\eta_1(t),\eta_2(t))\, \dot{y}(t)^2,\nonumber
\end{align}
where ${g_{i,j}:V\subseteq \R^2 \to \R}$, $V = [0, I_2] \Phi(U_\beta)$ open, are appropriate functions with $g_{1,0}(0,0) = g_{2,0}(0,0)=0$. The explicit representation can be found in Appendix~\ref{Sec:Appendix}. Here we highlight that the internal dynamics depend nonlinearly on~$\dot{y}$.

\subsection{Controller design}

The controller design is inspired by~\cite{Berg20a}, and we apply these results to the internal dynamics. To this end, we linearize the internal dynamics \eqref{eq:ID} around the equilibrium~${(\eta_1,\eta_2) =(0,0)}$, $\dot{y}=0$ and obtain
\begin{align}
\label{eq:ID-Linearization}
\begin{pmatrix} \dot{\eta}_1(t) \\ \dot{\eta}_2(t) \end{pmatrix} &= Q \begin{pmatrix} \eta_1(t) \\ \eta_2(t)  \end{pmatrix} + P \dot{y}(t),
\end{align}
where
\begin{align*}
Q = \begin{bmatrix} 0 & -12 \\ \frac{-c}{l^2m} & \frac{12d}{l^2m} \end{bmatrix}, \ P = \begin{bmatrix} 10 \\ \frac{-10d}{l^2m} \end{bmatrix}
\end{align*}
and the matrix~$Q$ has eigenvalues
\begin{align*}
\lambda_1 &= \frac{6d}{l^2m} - 2\sqrt{\left(\frac{3d}{l^2m}\right)^2 + \frac{3c}{l^2m}}, \\
\lambda_2 &= \frac{6d}{l^2m} + 2\sqrt{\left(\frac{3d}{l^2m}\right)^2 + \frac{3c}{l^2m}}.
\end{align*}
Note, that for ${c > 0}$ we have ${\lambda_1 < 0 < \lambda_2}$, thus~\eqref{eq:ID-Linearization} has a hyperbolic equilibrium. Therefore, the linearized internal dynamics have an unstable part.
Now, we find a transformation
\begin{align*}
V = \begin{bmatrix} \tfrac{\lambda_1 l^2 m}{c} & \tfrac{\lambda_2 l^2 m}{c} \\ 1 & 1 \end{bmatrix}\in \Gl_2(\R),
\end{align*}
which diagonalizes~$Q$ (and hence separates the stable and the unstable part of the internal dynamics) such that
\begin{align*}
V^{-1} Q V = \begin{bmatrix} \lambda_1 & 0 \\ 0 & \lambda_2 \end{bmatrix}, \quad V^{-1}P = \begin{bmatrix} p_1 \\ p_2 \end{bmatrix},
\end{align*}
where $p_{1,2} = \pm \tfrac{10}{Dc} (c + d \lambda_{1,2})$, with~${D:= \frac{l^2 m}{c} (\lambda_1 - \lambda_2)} = \det(V)$.
Using the transformation $\hat{\eta}(t) = V^{-1} \eta(t)$ we obtain the linearized internal dynamics separated in a stable and an unstable part
\begin{equation}\label{eq:ID-separated}
\begin{aligned}
\dot{\hat{\eta}}_1(t) &= \lambda_1 \hat{\eta}_1(t) +\tfrac{10}{Dc} (c + d \lambda_1) \dot y(t),\\
\dot{\hat{\eta}}_2(t) &= \lambda_2 \hat{\eta}_2(t) - \tfrac{10}{Dc} (c + d \lambda_2) \dot{y}(t).
\end{aligned}
\end{equation}
In virtue of~\cite{Berg20a} we seek to define an auxiliary output $\hat{y}_{\rm new}$ as a flat output for system~\eqref{eq:ID-separated}, i.e., its second equation, and calculate the new relative degree of system~\eqref{eq:motion-manipulator} with respect to $\hat{y}_{\rm new}$.
However, it is not obvious how to treat this task. One way would be to express the new output in terms of the original coordinates (i.e., replace~$\eta_2$ in~\eqref{eq:ID-Linearization} with $\eta_2$ in~\eqref{eq:ID} and perform the same transformations which lead to~$\hat\eta_2$) and calculate the derivatives. Another way, and this is the one we choose, is to use the linearization of the internal dynamics to calculate the relative degree based on~\eqref{eq:ID-separated}. To avoid confusion we do not use the dot for time derivatives here, but a superscript to indicate the derivatives w.r.t.\ the linearization~\eqref{eq:ID-separated}. The new output $\hat{y}_{\rm new}$ may be given by the variable of the unstable part as~${\hat{y}_{\rm new} = \hat{\eta}_2}$, and we calculate
\begin{align*}
\hat{y}^{[1]}_{\rm new}(t) &= \lambda_2 \hat{\eta}_2(t) + p_2 \dot{y}(t), \\
\hat{y}^{[2]}_{\rm new}(t) &= \lambda_2 ( \hat{\eta}_2(t) +  p_2 \dot{y}(t) ) + p_2 \ddot{y}(t).
\end{align*}
In the last equation we may insert the equation for~$\ddot y$ according to~\eqref{eq:motion-manipulator}, \eqref{eq:output-manipulator}, which explicitly depends on the input~$u$. Hence, the relative degree of~\eqref{eq:motion-manipulator} with respect to~$\hat{y}_{\rm new}$ is again $r_{\rm new}=2$, thus remains unchanged. Recalling the aim formulated in Section~\ref{Ssec:Methodology}, we need to find a different output. To this end, we use the transformation~${\bar{\eta}_2(t) = \hat{\eta}_2(t) - p_2 y(t)}$ to obtain
\begin{align}
\label{eq:DGL-eta2}
\dot{\bar{\eta}}_2(t) = \lambda_2 \bar{\eta}_2(t) + \lambda_2 p_2 y(t)
\end{align}
and define
\begin{align}
\label{eq:y-new}
y_{\rm new}(t) = \bar{\eta}_2(t). 
\end{align}
A short calculation shows that system~\eqref{eq:motion-manipulator} with new output~\eqref{eq:y-new} has relative degree~$r_{\rm new}=3$:
\begin{align}
\label{eq:dot_ynew-approx}
y^{[1]}_{\rm new}(t) &\overset{\eqref{eq:DGL-eta2}}{=} \lambda_2 (\bar{\eta}_2(t) + p_2 y(t)) \nonumber \\
y^{[2]}_{\rm new}(t) &= \lambda_2^2( \bar{\eta}_2(t) +  p_2 y(t)) + \lambda_2 p_2 \dot{y}(t) \\
y^{[3]}_{\rm new}(t) &= \lambda_2^3( \bar{\eta}_2(t) +  p_2 y(t)) + \lambda_2^2 p_2 \dot{y}(t) + \lambda_2 p_2 \ddot{y}(t), \nonumber
\end{align}
where in the last equation we again use~\eqref{eq:motion-manipulator},~\eqref{eq:output-manipulator} to obtain the input~$u$ explicitly.

\begin{Rem}
We stress that the derivatives~$y_{\rm new}^{[i]}$  are calculated by replacing the original unstable part of the internal dynamics, i.e., the second equation in~\eqref{eq:ID}, by the linearized version~\eqref{eq:DGL-eta2}. They are not equivalent to those expressed in the original coordinates, which means that ${y^{[1]}_{\rm new}(t) \neq \dot{y}_{\rm new}(t)}$. More precisely, using the original coordinates would mean that the transformations leading from~$\eta_2$ to~$\bar \eta_2$ are performed with~$\eta_2$ in~\eqref{eq:ID} so that
\begin{equation}\label{eq:bareta2-orig}
    \bar \eta_2 = \hat \eta_2 - p_2 y = [0,1] V^{-1} \begin{pmatrix}\eta_1\\ \eta_2\end{pmatrix} - p_2 y,
\end{equation}
and $\eta_1, \eta_2, y$ are replaced via~\eqref{eq:BIF-Manipulator}.
\end{Rem}

Now, in order to track the original reference with the original output, we have to find a new reference signal for system~\eqref{eq:motion-manipulator} with new output~\eqref{eq:y-new}.
The new reference~${\bar{y}_{\rm ref}}$ is given by the solution of~\eqref{eq:DGL-eta2} when the original output~$y$ is substituted by the original reference signal~$y_{\rm ref}$:
\begin{align}
\label{eq:ref-new}
\dot{\bar{\eta}}_{2,\rm ref}(t) &= \lambda_2 \bar{\eta}_{2,\rm ref}(t) + \lambda_2 p_2 y_{\rm ref}(t), &\bar{\eta}_{2,\rm ref}(0) = \bar{\eta}_{2,\rm ref}^0 \nonumber \\
\bar{y}_{\rm ref}(t) &= \bar{\eta}_{2,\rm ref}(t).
\end{align}
We stress that~\eqref{eq:ref-new} adds a dynamic equation to the overall controller design. In order for the controller from~\cite{BergLe18a} to be applicable,~$\bar{y}_{\rm ref}$ and its derivatives should be bounded. To show this, we use the following well known result: there exists ${x \in W^{1,\infty}(\R_{\geq 0} \to \R)}$ solving ${\dot{x}(t) = \lambda x(t) + \gamma g(t)}$, where $\lambda > 0, \gamma \in \R$, ${g \in L^{\infty}(\R_{\geq 0} \to \R)}$, and~${x(0)=x^0}$ if, and only if, $x^0 + \int_0^\infty e^{-\lambda s} \gamma g(s) ds = 0$. Hence we set
\begin{align}
\label{eq:eta2-ref-0}
\bar{\eta}_{2,\rm ref}^0 = -\int_0^\infty e^{-\lambda_2 s} \lambda_2 p_2 y_{\rm ref}(s) ds.
\end{align}
We emphasize that, if~$y_{\rm ref}$ is generated by an exosystem
\[
    \dot w(t) = A_e w(t),\quad  y_{\rm ref}(t) = C_e w(t),\quad w(0) = w^0,
\]
with known parameters~$A_e\in\R^{k\times k}, C_e\in\R^{1\times k}$ and $w^0\in\R^k$ such that $\sigma(A_e)\subseteq \overline{\C_-}$ and any eigenvalue $\lambda\in\sigma(A_e)\cap i\R$ is semisimple (note that this guarantees $y_{\rm ref}\in W^{1,\infty}(\R_{\ge 0}\to\R)$), then $\bar{\eta}_{2,\rm ref}^0$ can be calculated via the solution~$X\in\R^{1\times k}$ of the Sylvester equation
\[
    Q_2 X - X A_e = P_2 C_e,
\]
where in our example $Q_2=\lambda_2$ and $P_2 = \lambda_2 p_2$. It is shown in~\cite[Lem.~3.2]{Berg20a} that in this case
\[
    \bar{\eta}_{2,\rm ref}^0 = -Xw^0.
\]
The tracking error for the auxiliary system with new output~\eqref{eq:y-new} and new reference~\eqref{eq:ref-new} is defined by
\begin{align}
\label{eq:e-new}
e_0(t) = y_{\rm new}(t) - \bar{y}_{\rm ref}(t).
\end{align}
Applying the funnel control law from~\cite{BergLe18a} requires the derivatives of~$y_{\rm new}$, and in order to implement it we have to express them in terms of the original variables~$q$. In the following we present two approaches to obtain these derivatives. In the first approach we use the representation in~\eqref{eq:dot_ynew-approx} and replace~$\bar\eta_2$ by the original variables via~\eqref{eq:bareta2-orig} and~\eqref{eq:BIF-Manipulator}. In the second approach we only replace~$y_{\rm new}$ by the original variables and approximate the derivatives of this representation using a high-gain observer.

For later use and $\varphi_0\in\Phi_3$, $\kappa_0>0$, we define the expressions
\begin{align}
\label{notation-derivative-e}
e^{[1]}_0 &= y^{[1]}_{\rm new} - \dot{\bar{y}}_{\rm ref}, \nonumber\\
e^{[2]}_0 &= y^{[2]}_{\rm new} - \ddot{\bar{y}}_{\rm ref}, \nonumber\\
k^{[1]}_0 &= \frac{2 \kappa_0 \varphi_0 e_0}{(1 - \varphi_0^2 e_0^2)^2} \left(\dot{\varphi}_0 e_0 + \varphi_0 e^{[1]}_0 \right), \\
e^{[1]}_1 &= e^{[2]}_0 + k_0 e^{[1]}_0 + k^{[1]}_0 e_0. \nonumber
\end{align}

\subsubsection{Linearization.}
The first option is to replace~$\bar\eta_2$ in~$y_{\rm new}$ and its derivatives in~\eqref{eq:dot_ynew-approx} by~\eqref{eq:bareta2-orig} and to express everything in terms of the original coordinates via~~\eqref{eq:BIF-Manipulator}. To this end, we consider
\begin{align*}
\Psi: U_\beta &\to \R \\
 (x_1,\ldots,x_4)
&\mapsto
- p_2 (x_1 + \tfrac{1}{2} x_2) \\
& + \tfrac{1}{D} \left( -x_2 + \tfrac{\lambda_2 l^2 m}{c} \left[ (\tfrac{1}{3} + \tfrac{1}{2} \cos(x_2))x_3 + \tfrac{1}{3} x_4 \right] \right)
\end{align*}
and replace $\bar\eta_2(t) = \Psi(\alpha(t),\beta(t),\dot\alpha(t),\dot\beta(t))$ in~$y_{\rm new}(t) = \bar\eta_2(t)$ and in~\eqref{eq:dot_ynew-approx}, thus
\begin{align}
\label{eq:derivatives-ynew-approx}
y_{\rm new}(t) &= \Psi(\alpha(t),\beta(t),\dot\alpha(t),\dot\beta(t)),\nonumber \\
y_{\rm new}^{[1]}(t) &= \lambda_2 \Psi(\alpha(t),\beta(t),\dot\alpha(t),\dot\beta(t)) + \lambda_2 p_2 \big(\alpha(t) + \tfrac{1}{2} \beta(t)\big),\nonumber \\
y_{\rm new}^{[2]}(t) &= \lambda_2^2 \Psi(\alpha(t),\beta(t),\dot\alpha(t),\dot\beta(t)) + \lambda_2^2 p_2 \big(\alpha(t) + \tfrac{1}{2} \beta(t)\big)  \nonumber\\
&\quad + \lambda_2 p_2 \big(\dot \alpha(t) + \tfrac{1}{2} \dot\beta(t)\big).
\end{align}
With this, the application of the controller from~\cite{BergLe18a} to system~\eqref{eq:motion-manipulator} with new output~\eqref{eq:y-new} and reference signal as in~\eqref{eq:ref-new} leads to the following overall control law:

\fbox{\parbox{0.46\textwidth}{%
\begin{align}
\label{eq:Control}
\dot{\bar{\eta}}_{2,\rm ref}(t) &= \lambda_2 \bar{\eta}_{2,\rm ref}(t) + \lambda_2 p_2 y_{\rm ref}(t), \quad \bar{\eta}_{2,\rm ref}(0) = \bar{\eta}_{2,\rm ref}^0, \nonumber \\
\bar{y}_{\rm ref}(t) &= \bar{\eta}_{2,\rm ref}(t), \nonumber \\
y_{\rm new}(t) &= \Psi(\alpha(t),\beta(t),\dot\alpha(t),\dot\beta(t)),\nonumber \\
e_0(t) &= y_{\rm new}(t) - \bar{y}_{\rm ref}(t), \nonumber \\ 
e_1(t) &= e^{[1]}_0(t) + k_0(t)e_0(t),\quad e^{[1]}_0\ \text{via~\eqref{notation-derivative-e} and~\eqref{eq:derivatives-ynew-approx}}, \nonumber \\
e_2(t) &= e^{[1]}_1(t) + k_1(t)e_1(t),\quad e^{[1]}_1\ \text{via~\eqref{notation-derivative-e} and~\eqref{eq:derivatives-ynew-approx}}, \nonumber \\
k_i(t) &= \frac{1}{1 - \varphi_i(t)^2 e_i(t)^2},\quad i=0,1,2, \nonumber \\
u(t) &= k_2(t) e_2(t),
\end{align}}}
with~$\bar{\eta}_{2,\rm ref}^0$ as in~\eqref{eq:eta2-ref-0}, reference signal $y_{\rm ref}\in W^{1,\infty}(\R_{\ge 0}\to\R)$ and funnel functions $\varphi_i \in\Phi_{3-i}$ for $i=0,1,2$.
Note that, since~$\Gamma(x)$ is negative for $x \in U_\beta$, the control input~$u$ has a positive sign according to~\cite{BergLe18a}.

\subsubsection{High-gain observer.}
The second option is to approximate the first and second derivative of $y_{\rm new}$ using a high-gain observer, see~\cite{Khal01}. To this end, consider
\begin{align*}
\dot{\zeta}(t) =  L (y_{\rm new}(t) - \zeta_1(t)) + Z \zeta(t) , \quad \zeta(0) = \zeta^0 \in \R^{3},
\end{align*}
where
\begin{equation*}
Z = \begin{bmatrix}0 & I_{2} \\ 0 &0 \end{bmatrix}, \quad L = \begin{bmatrix} l_1 \\ l_2 \\ l_3 \end{bmatrix},
\end{equation*}
with~${l_i \in \R}$ for~$i=1,2,3$. Then $\zeta_2$ and $\zeta_3$ approximate the first and second derivative of $y_{\rm new}$, thus for the controller we set $y_{\rm new}^{[1]} = \zeta_2$ and $y_{\rm new}^{[2]} = \zeta_3$. The high-gain observer is a dynamical part of the controller and invoking~\eqref{notation-derivative-e} the overall controller reads:

\fbox{\parbox{0.46\textwidth}{%
\begin{align}
\label{eq:Control-HG}
\dot{\bar{\eta}}_{2,\rm ref}(t) &= \lambda_2 \bar{\eta}_{2,\rm ref}(t) + \lambda_2 p_2 y_{\rm ref}(t), && \bar{\eta}_{2,\rm ref}(0) = \bar{\eta}_{2,\rm ref}^0, \nonumber \\
\bar{y}_{\rm ref}(t) &= \bar{\eta}_{2,\rm ref}(t), \nonumber \\
y_{\rm new}(t) &= \Psi(\alpha(t),\beta(t),\dot\alpha(t),\dot\beta(t)),\nonumber \\
\dot{\zeta}_1(t) &= l_1 (y_{\rm new}(t) - \zeta_1(t)) + \zeta_2(t) ,&& \zeta_1(0) = \zeta_1^0\nonumber \\
\dot{\zeta}_2(t) &= l_2 (y_{\rm new}(t) - \zeta_1(t)) + \zeta_3(t) ,&&  \zeta_2(0) = \zeta_2^0 \nonumber \\
\dot{\zeta}_3(t) &= l_3 (y_{\rm new}(t) - \zeta_1(t)), &&  \zeta_3(0) = \zeta_3^0 \nonumber \\
y^{[1]}_{\rm new}(t) &= \zeta_2(t), \nonumber \\
y^{[2]}_{\rm new}(t) &= \zeta_3(t), \nonumber \\
e_0(t) &= y_{\rm new}(t) - \bar{y}_{\rm ref}(t) \nonumber \\
e_1(t) &= e^{[1]}_0(t) + k_0(t)e_0(t),&& e^{[1]}_0\ \text{via~\eqref{notation-derivative-e}}, \nonumber \\
e_2(t) &= e^{[1]}_1(t) + k_1(t)e_1(t),&& e^{[1]}_1\ \text{via~\eqref{notation-derivative-e}}, \nonumber \\
k_i(t) &= \frac{1}{1 - \varphi_i(t)^2 e_i(t)^2}, && i=0,1,2, \nonumber \\
u(t) &= k_2(t) e_2(t),
\end{align}
}}

with~$\bar{\eta}_{2,\rm ref}^0$ as in~\eqref{eq:eta2-ref-0}, high-gain observer parameters~$l_i \in \R$ and initial values $\zeta_i^0\in\R$, $i=1,2,3$, reference signal $y_{\rm ref}\in W^{1,\infty}(\R_{\ge 0}\to\R)$ and funnel functions $\varphi_i \in\Phi_{3-i}$ for $i=0,1,2$.

\section{Numerical case study}
\label{NumericalExample}
In this section we present the results of the numerical case study.
We perform end-effector output tracking of the manipulator arm~\eqref{eq:motion-manipulator} with $l=1\, \mathrm{m}$, $m=1\, \mathrm{kg}$, $c=1\, \mathrm{Nm/rad}$ and $d=0.25\, \mathrm{Nm s/rad}$ applying controls~\eqref{eq:Control} and~\eqref{eq:Control-HG}, resp.
As a reference signal we choose the trajectory from~\cite{SeifBlaj13}
\begin{equation*} \label{eq:reference}
\begin{aligned}
y_{\rm ref}(t) = y_0 &+ \left[ 126 (\tfrac{t}{t_f - t_0})^5 - 420(\tfrac{t}{t_f - t_0})^6 + 540(\tfrac{t}{t_f - t_0})^7\right. \\
 &\quad\ \ \left.- 315(\tfrac{t}{t_f - t_0})^8 + 70(\tfrac{t}{t_f - t_0})^9 \right](y_f - y_0),
\end{aligned}
\end{equation*}
which establishes a transition from~$y_0$ to $y_f$ within the time interval~$t_0$ to~$t_f$.
%
The situation is depicted in Figure~\ref{Fig:Transition}.
\begin{center}
\centering
\resizebox{1\linewidth}{!}{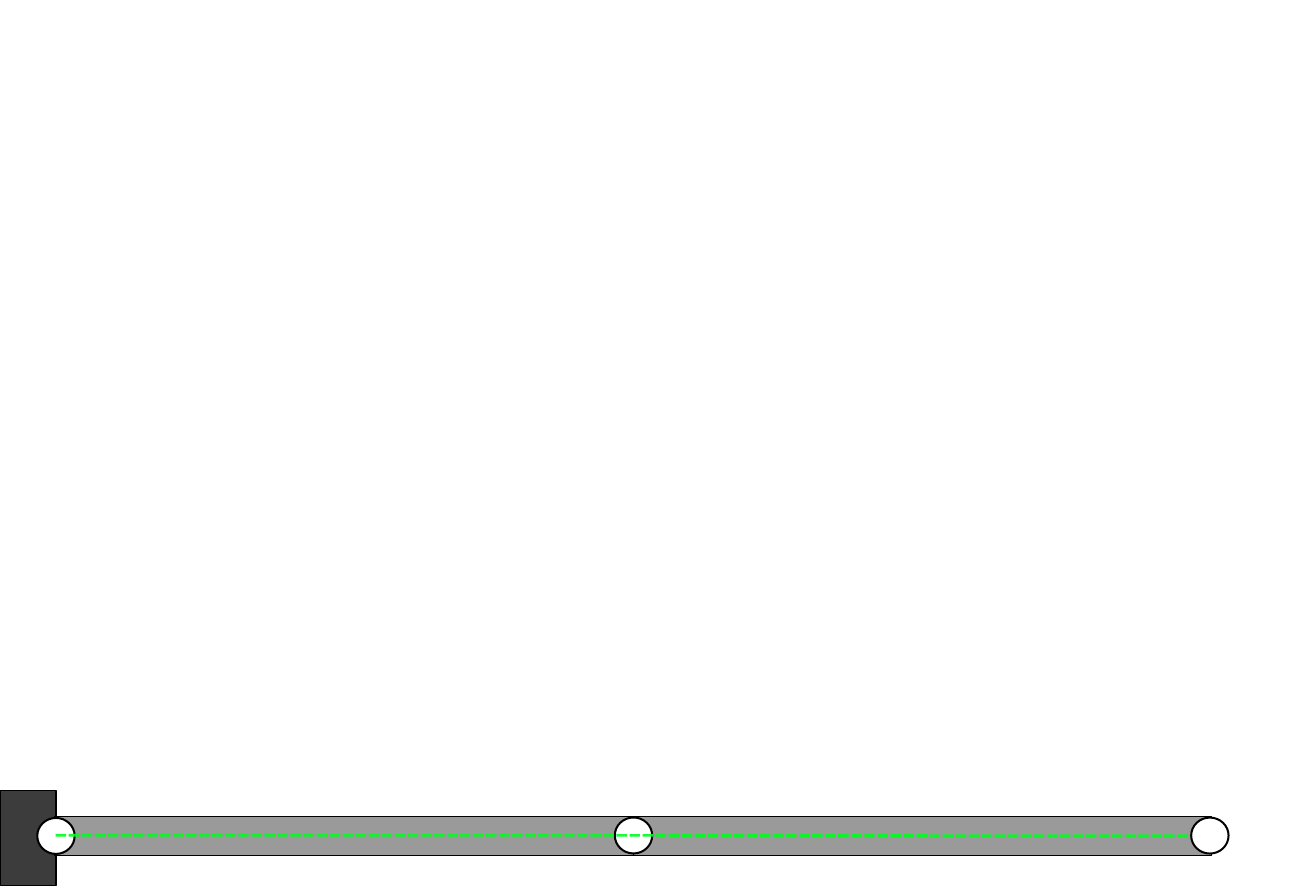}
\captionof{figure}{Schematic snapshots of the manipulator arm's transition established by the reference trajectory $y_{\rm ref}$ with $y_0 = 0 \, \rm rad$ and $y_f=\tfrac{\pi}{4} \, \rm rad$.}
\label{Fig:Transition}
\end{center}
We choose $y_0 = 0\, \mathrm{rad}$,
$y_f = \tfrac{\pi}{4}\, \mathrm{rad}$,
 $t_0 = 0 \,\mathrm{s}$ and
$t_f = 3\, \mathrm{s}$, which means the transition is performed in rather short time, namely within three seconds.
%
We choose the funnel functions
\begin{align*}
&\varphi_0(t) = \varphi_1(t) = \left(1.5 \, e^{-0.8 \, t} + 0.001 \right)^{-1}, \\
&\varphi_2(t) = \left(60 e^{-0.2 \, t} + 0.001\right)^{-1}, \quad t \ge 0,
\end{align*}
which ensure that the initial errors $e_0(0), e_1(0), e_2(0)$ lie within the respective funnel boundaries;
and we choose the high-gain parameters $l_1 = 10^2, l_2=10^5, l_3 = 10^6$.
Furthermore, we assume the control torque to be affected by high frequent disturbances~$d(t)$ as expectable in real applications caused e.g. by
friction or unexpected external vibrating forces, cf.~\cite[Ch.~11~\&~13]{Hack17}.
For simulation purposes we choose $d(t) = 0.1 \sin(5t) + 0.2\cos(8t)$ and recall that in~\eqref{eq:motion-manipulator} we have $u_d(t) = u(t) + d(t)$.
%
%
\begin{figure}[htbp]
\subfloat[Error~$e_0$ and funnel boundary~$\varphi_0^{-1}$.]{%
\includegraphics[width=0.95\linewidth]{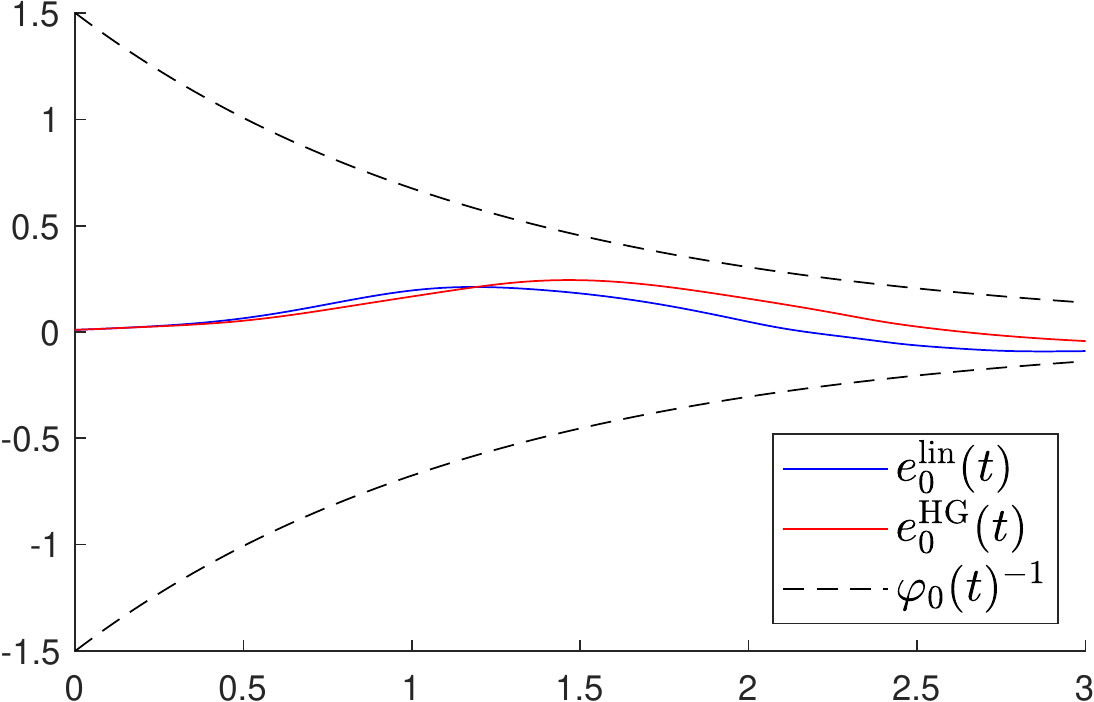}
\label{Fig:e0}
} \\
\subfloat[Angles $\alpha$ and $\beta$ of the manipulator arm.]{%
\includegraphics[width=0.95\linewidth]{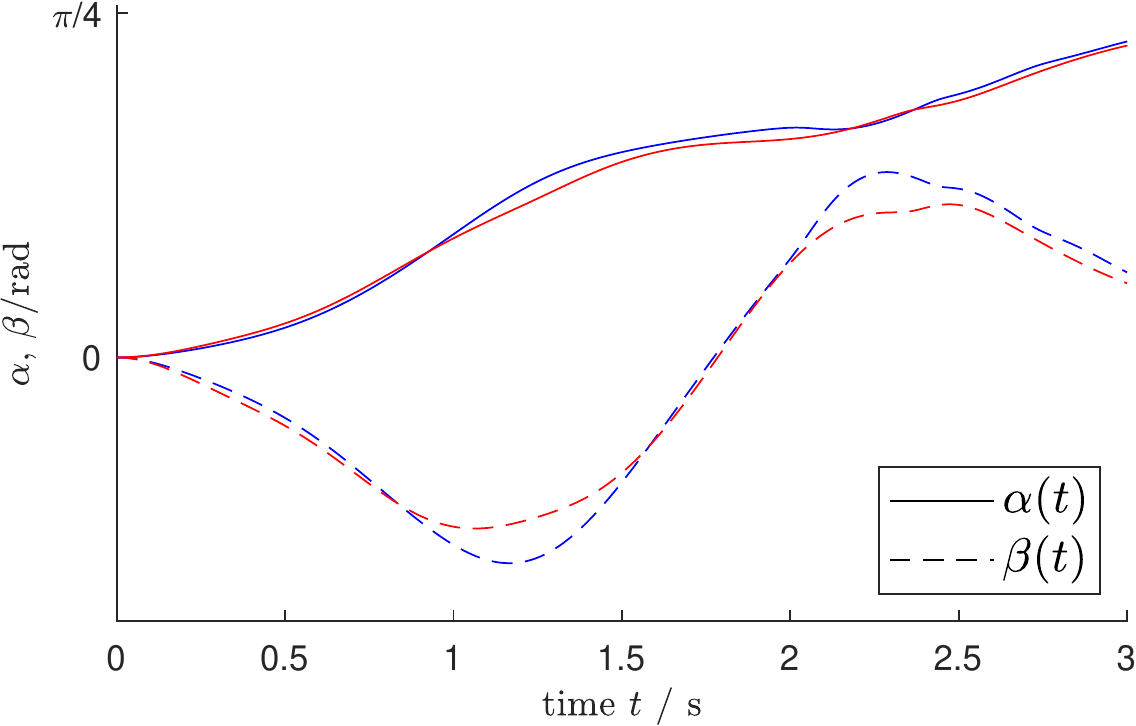}
\label{Fig:states}
} \\
\subfloat[Input function $u$.]{%
\includegraphics[width=0.95\linewidth]{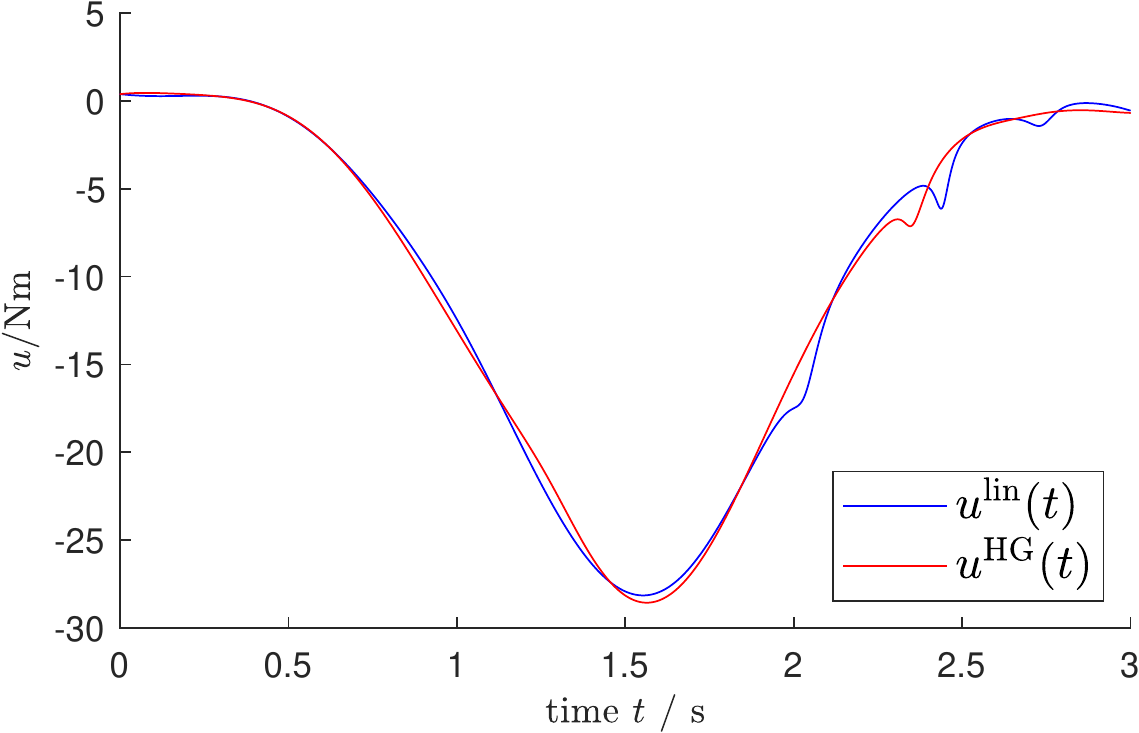}
\label{Fig:control}
}
\caption{Simulation of the controllers~\eqref{eq:Control} (superscript~\textit{lin}) and~\eqref{eq:Control-HG} (superscript~\textit{HG}) applied to~\eqref{eq:motion-manipulator},~\eqref{eq:output-manipulator} under disturbances~$d$.}
\label{Fig:Simulation}
\end{figure}

Figure~\ref{Fig:Simulation} shows the results of the simulations, which have been performed in \textsc{Matlab}
(solver: \textsf{ode15s}, rel. tol.: $10^{-9}$, abs. tol.: $10^{-12}$)
over the time interval $0-3\, \mathrm{s}$.
In the following the superscript~\textit{lin} refers to the results under the control~\eqref{eq:Control}, where the derivatives of~$y_{\rm new}$ are expressed via~\eqref{eq:derivatives-ynew-approx}; the superscript~\textit{HG} refers to the results under the control~\eqref{eq:Control-HG}, where the derivatives of~$y_{\rm new}$ are approximated via a high-gain observer.
As depicted in Figure~\ref{Fig:e0} the error~$e_0$ evolves within the funnel boundary~$\varphi_0^{-1}$.
Figure~\ref{Fig:states} shows the state variables, i.e., the angles $(\alpha, \beta)$, of system~\eqref{eq:motion-manipulator} during the tracking process.
The large angle~$\beta$, i.e., large deflection of the passive link, results from the fast transition within three seconds.
However, note that ${\beta(t) \in \B}$ for $0 \leq t \leq 3$.
In Figure~\ref{Fig:control} the respective control inputs~$u$ are depicted.
Note that both approaches, the linearization and the use of a high-gain observer, generate comparable control inputs.
Figure~\ref{Fig:y_and_yref} shows that the goal of the numerical case study, namely end-effector output tracking of a prescribed trajectory, is successful.
\begin{center}
\includegraphics[width=0.95\linewidth]{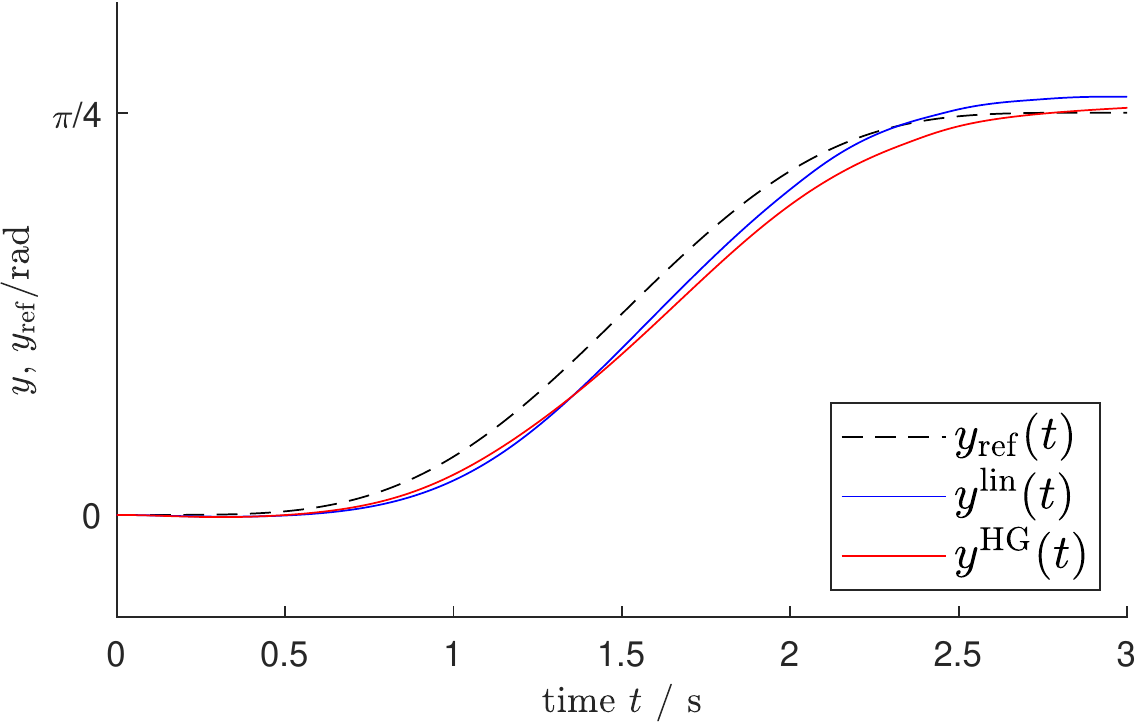}
\captionof{figure}{Output~$y$ and reference~$y_{\rm ref}$.}
\label{Fig:y_and_yref}
\end{center}

\section{Conclusion}
In the present paper we exploited ideas and methods recently found for linear non-minimum phase systems in~\cite{Berg20a} to proceed a numerical case study for a nonlinear non-minimum phase robotic manipulator. To this end, we transformed the system under consideration into Byrnes-Isidori form in order to decouple the internal dynamics. We linearized the internal dynamics around an equilibrium point and separated its stable and unstable parts. Then we followed the steps in the aforementioned work in order to design the feedback control laws~\eqref{eq:Control} and \eqref{eq:Control-HG} using different methods to approximate the new output's derivatives. Utilizing these control laws we simulated end-effector output tracking of the rotational manipulator arm~\eqref{eq:motion-manipulator}. The simulation showed that output tracking of such a nonlinear non-minimum phase system is successful, even under the action of disturbances.

We stress that the present paper is only a case study in order to gain some insight as to whether the methods and ideas derived in~\cite{Berg20a} can be extended to nonlinear non-minimum phase systems. The presented results give reason to adopt the presented techniques to achieve efficient control methods for nonlinear systems with unstable internal dynamics via funnel control.


\appendix

\section{The functions \MakeLowercase{$g_{i,j}$}}\label{Sec:Appendix}

Here we list the functions $g_{i,j}$ that appear in~\eqref{eq:ID}.
\begin{align*}
g_{1,0} : V \subseteq \R^2 &\to \R, \\
(\eta_1,\eta_2) &\mapsto \frac{12}{2-3 \cos(\eta_1)} \eta_2, \\
\ \\
g_{1,1} : V \subseteq \R^2 &\to \R, \\
(\eta_1,\eta_2) &\mapsto \frac{-12}{2 - 3\cos(\eta_1)}\left(\frac{1}{3} + \frac{1}{2} \cos(\eta_1) \right), \\
\ \\
g_{2,0} : V \subseteq \R^2 &\to \R, \\
(\eta_1,\eta_2) &\mapsto \frac{-1}{l^2 m (2 - \cos(\eta_1))^2} \\
&\cdot \biggl\{
\eta_1 \big[4c - 12 c + 9 c \cos(\eta_1)^2 \big] \\
&+ 4\eta_2 \big[6 d - 9 \big( d \cos(\eta_1) + l^2 m \eta_2 \sin(\eta_1) \big) \big] \biggr\}, \\
\ \\
g_{2,1} : V \subseteq \R^2 &\to \R, \\
(\eta_1,\eta_2) &\mapsto \frac{2}{l^2 m (2 - \cos(\eta_1))^2 (16 - 9 \cos(\eta_1)^2)} \\
&\cdot \biggl\{ \eta_2 \big[ 453 \sin(\eta_1) + 216 \sin(2 \eta_1) - 27 \sin(3 \eta_1) \\
&\quad+  l^2 m  \big( \tfrac{-1071}{2} \sin(\eta_1) - \tfrac{1071}{4} \sin(2 \eta_1) \\
&+ \tfrac{81}{4} \sin(3 \eta_1) + \tfrac{81}{8} \sin(4 \eta_1) \big) \big] \\
&+  \tfrac{35d}{8} - \tfrac{99 d }{2} \cos(2 \eta_1) + \tfrac{81 d}{8} \cos(4 \eta_1) \biggr\}, \\
\ \\
g_{2,2}: V \subseteq \R^2 &\to \R, \\
(\eta_1,\eta_2) &\mapsto \frac{-4 \sin(\eta_1)}{l^2 m (2 - \cos(\eta_1))^2 (16 - 9 \cos(\eta_1)^2)} \\
&\cdot \biggl\{
80 + 144\cos(\eta_1) + 90 \cos(\eta_1) \\
&- \big[80 + 192 \cos(\eta_1) + 90 \cos(\eta_1)^2 \\
&- 27 \cos(\eta_1)^3 \big] l^2 m \biggr\}.
\end{align*}

\end{document}

%% file: RoboterArm.pdf_tex
\begingroup%
  \makeatletter%
  \providecommand\color[2][]{%
    \errmessage{(Inkscape) Color is used for the text in Inkscape, but the package 'color.sty' is not loaded}%
    \renewcommand\color[2][]{}%
  }%
  \providecommand\transparent[1]{%
    \errmessage{(Inkscape) Transparency is used (non-zero) for the text in Inkscape, but the package 'transparent.sty' is not loaded}%
    \renewcommand\transparent[1]{}%
  }%
  \providecommand\rotatebox[2]{#2}%
  \newcommand*\fsize{\dimexpr\f@size pt\relax}%
  \newcommand*\lineheight[1]{\fontsize{\fsize}{#1\fsize}\selectfont}%
  \ifx\svgwidth\undefined%
    \setlength{\unitlength}{378.49636898bp}%
    \ifx\svgscale\undefined%
      \relax%
    \else%
      \setlength{\unitlength}{\unitlength * \real{\svgscale}}%
    \fi%
  \else%
    \setlength{\unitlength}{\svgwidth}%
  \fi%
  \global\let\svgwidth\undefined%
  \global\let\svgscale\undefined%
  \makeatother%
  \begin{picture}(1,0.67519171)%
    \lineheight{1}%
    \setlength\tabcolsep{0pt}%
    \put(0,0){\includegraphics[width=\unitlength,page=1]{RoboterArm.pdf}}%
    \put(0.9137644,0.03307029){\color[rgb]{0,0,0}\makebox(0,0)[lt]{\lineheight{1.25}\smash{\begin{tabular}[t]{l}S\end{tabular}}}}%
    \put(0,0){\includegraphics[width=\unitlength,page=2]{RoboterArm.pdf}}%
    \put(0.6565363,0.65271839){\color[rgb]{0,0,0}\makebox(0,0)[lt]{\lineheight{1.25}\smash{\begin{tabular}[t]{l}S\end{tabular}}}}%
    \put(0,0){\includegraphics[width=\unitlength,page=3]{RoboterArm.pdf}}%
    \put(0.94340332,0.02983189){\color[rgb]{0,0,0}\makebox(0,0)[lt]{\lineheight{1.25}\smash{\begin{tabular}[t]{l}$t=t_0$\end{tabular}}}}%
    \put(0.68525669,0.65540783){\color[rgb]{0,0,0}\makebox(0,0)[lt]{\lineheight{1.25}\smash{\begin{tabular}[t]{l}$t=t_f$\end{tabular}}}}%
    \put(0.88824573,0.40568736){\color[rgb]{0,0,0}\makebox(0,0)[lt]{\lineheight{1.25}\smash{\begin{tabular}[t]{l}$y_{\rm ref}$\end{tabular}}}}%
    \put(0,0){\includegraphics[width=\unitlength,page=4]{RoboterArm.pdf}}%
    \put(0.90104473,0.05654012){\color[rgb]{0,0,0}\makebox(0,0)[lt]{\lineheight{1.25}\smash{\begin{tabular}[t]{l}S\end{tabular}}}}%
    \put(0.75819928,0.50859451){\color[rgb]{0,0,0}\makebox(0,0)[lt]{\lineheight{1.25}\smash{\begin{tabular}[t]{l}S\end{tabular}}}}%
  \end{picture}%
\endgroup%